\documentclass{stavanger}

\usepackage[utf8]{inputenc}
\usepackage{amssymb,amsmath,amsfonts,amsthm}
\usepackage{ dsfont }
\usepackage{hyperref}
\usepackage{natbib}
\usepackage{graphicx}
\usepackage{tikz}
\usepackage{cleveref}
\usepackage{mathtools}
\usepackage{placeins}
\usepackage{bm}
\usepackage{eufrak}
\usepackage{animate}
\usepackage[export]{adjustbox} 
\usepackage{mathbbol}

\startlocaldefs

\endlocaldefs

\begin{document}


\begin{frontmatter}

\begin{fmbox}


\dochead{Preprint}{FP}



\title{The crucial role of Lagrange multipliers in a space-time symmetry preserving discretization scheme for IVPs}


\author[
   addressref={aff1},                   	  
   corref={aff1},                     		  
   email={alexander.rothkopf@uis.no}   		  
]{\inits{AR}\fnm{Alexander} \snm{Rothkopf}}
\author[
   addressref={aff2,aff3},                   	  
   email={jan.nordstrom@liu.se}   		  
]{\inits{JN}\fnm{Jan} \snm{Nordstr{\"o}m}}


\address[id=aff1]{
  \orgname{Faculty of Science and Technology}, 	 
  \street{University of Stavanger},                     		 
  \postcode{4021},                               			 
  \city{Stavanger},                              				 
  \cny{Norway}                                   				 
}

\address[id=aff2]{
  \orgname{Department of Mathematics}, 	 
  \street{Link{\"o}ping University},                     		 
  \postcode{SE-581 83},                               			 
  \city{Link{\"o}ping},                              				 
  \cny{Sweden}                                   				 
}

\address[id=aff3]{
  \orgname{Department of Mathematics and Applied Mathematics}, 	 
  \street{University of Johannesburg},                     		 
  \postcode{P.O. Box 524, Auckland Park 2006},                           
  \city{Johannesburg},                              				 
  \cny{South Africa}                                   				 
}



\end{fmbox}


\begin{abstractbox}
\begin{abstract} 
In a recently developed variational discretization scheme for second order initial value problems \cite{Rothkopf:2023ljz}, it was shown that the Noether charge associated with time translation symmetry is exactly preserved in the interior of the simulated domain. The obtained solution also fulfils the naively discretized equations of motions inside the domain, except for the last two grid points. Here we provide an explanation for the deviations at the boundary as stemming from the Lagrange multipliers used to implement initial and connection conditions. We show explicitly that the Noether charge including the boundary corrections is exactly preserved at its continuum value over the whole simulation domain, including the boundary points.
\end{abstract}



\end{abstractbox}


\end{frontmatter}


\section{Introduction}

The general theory of relativity underlies the world-line formalism to describe the dynamics of a point particle. It considers the trajectory of the particle as a one-dimensional submanifold of the (d+1) dimensional spacetime with $d$ space and one time directions. The space-time coordinates of the particle ${\bf x}(\gamma)=(t(\gamma),x(\gamma))^T$ are considered as dependent variables, parametrised by the world-line parameter $\gamma\in[\gamma_i,\gamma_f]$. The trajectory of a free particle traveling from initial coordinates $x_i,t_i$ to final coordinates $x_f,t_f$ in a space-time characterized by a given metric tensor $g_{\mu\nu}$, can be obtained from the critical point of the following world-line action functional
\begin{align}
    {\cal E}_{\rm BVP}=\int_{\gamma_i}^{\gamma_f}d\gamma E_{\rm BVP}[t,\dot t,x,\dot x]=\int_{\gamma_i}^{\gamma_f}d\gamma \frac{1}{2}\Big( g_{00} \Big(\frac{dt}{d\gamma}\Big)^2 + g_{11}\Big(\frac{dx}{d\gamma}\Big)^2 \Big)\label{eq:actnri}.
\end{align}
Here we have restricted ourselves to a $(1+1)d$ spacetime with time $x^0=t(\gamma)$ and $x^1=x(\gamma)$ coordinates and metric $g={\rm diag}[g_{00},g_{11}]={\rm diag}[c^2,-1]$. There are four independent degrees of freedom $t, \dot t=dt/d\gamma, x, \dot x=dx/d\gamma$ that appear in this reparameterization invariant action. Note that since $t$ and $x$ appear within derivatives, the action is manifestly invariant under constant space-time translations ${\bf x} \to {\bf x} + {\bf c}$.

In case that the particle is subject to an external potential $V(x)$ that is of less magnitude as the relativistic rest-energy of the particle $V(x)/mc^2\ll1$ ( a condition fulfilled in conventional non-relativistic initial value problems ) we can absorb the effect of the potential into a redefinition of the temporal component of the metric as $g_{00}=c^2+2V(x)/m$ (for more details of the construction of the continuum formalism see \cite{Rothkopf:2023ljz}).

In order to formulate the corresponding variational principle in the continuum world-line formalism, we deploy the doubled degrees of freedom formalism, discussed in \cite{galley_classical_2013} and put to practice in \cite{Rothkopf:2022zfb}. We introduce the forward $t_1,x_1$ and backward $t_2,x_2$ branches and define the following system action
\begin{align}
    \nonumber {\cal E}_{\rm IVP}&=\int_{\gamma_i}^{\gamma_f} d\gamma \, E_{\rm IVP}[t_1,\dot t_1,x_1,\dot x_1,t_2,\dot t_2,x_2,\dot x_2] \\
    &=\int_{\gamma_i}^{\gamma_f} d\gamma \Big\{ E_{\rm BVP}[t_1,\dot t_1,x_1,\dot x_1] - E_{\rm BVP}[t_2,\dot t_2,x_2,\dot x_2] \Big\}\label{eq:EIVPcont}.
\end{align}
We must make sure that the critical point of this action is equivalent to the solution of the corresponding Euler-Lagrange equations. This requires us to cancel several boundary terms arising during the variation of ${\cal E}_{\rm IVP}$. This is achieved by implementing not only fixed \textit{initial conditions} $x_{1,2}(0)=x_i$, $t_{1,2}(0)=t_i$, $\dot t_1(0)=\dot t_i$ and $\dot x_1(0)=\dot x_i$ but also the following \textit{connecting conditions} at the final $\gamma_f$:  $x_2(\gamma_f)=x_1(\gamma_f)$, $t_2(\gamma_f)=t_1(\gamma_f)$, $\dot x_2(\gamma_f)=\dot x_1(\gamma_f)$ and $ \dot t_2(\gamma_f)=\dot t_1(\gamma_f)$.

\section{Explicit treatment of initial and connecting conditions}
\label{sec:expllagrange}

The functional \cref{eq:EIVPcont} is manifestly time translation symmetric and encodes the dynamics of a point particle evolving in the presence of a non-relativistic potential as initial value problem, including implicit initial and connecting conditions. In preparation for the discretized formulation of our approach, we must include the initial and connecting conditions explicitly. To this end we add eight Lagrange multipliers to ${\cal E}_{\rm IVP}$, whose role is to explicitly implement the initial conditions $(\lambda_{1-4})$ and the connecting conditions at the end of the forward and backward branches of our doubled degree of freedom construction $(\lambda_{5-8})$ as shown below
\begin{align}
 \nonumber   {\cal E}_{\rm IVP}^{\lambda}&=\int_{\gamma_i}^{\gamma_f} d\gamma \Big\{ E_{\rm BVP}[t_1,\dot t_1,x_1,\dot x_1] - E_{\rm BVP}[t_2,\dot t_2,x_2,\dot x_2] \Big\} \\
    +& \nonumber\lambda_1\big({ t}_1(\gamma_i)-t_i\big)+\lambda_2\big(\dot { t}_1(\gamma_i)-\dot t_i\big)+\lambda_3\big({ x}_1(\gamma_i)-x_i\big) +\lambda_4\big(\dot{ x}_1(\gamma_i)-\dot x_i\big)\\
\nonumber +&\lambda_5\big({ t}_1(\gamma_f)-{ t}_2(\gamma_f)\big) + \lambda_6\big( \dot{ t}_1(\gamma_f)- \dot{ t}_2(\gamma_f)\big) \\
+&\lambda_7\big({ x}_1(\gamma_f)-{ x}_2(\gamma_f)\big)+\lambda_8\big( \dot{ x}_1(\gamma_f)- \dot{ x}_2(\gamma_f)\big) \label{eq:EIVPcontExpl}.
\end{align}

These additional terms that make reference to the dynamical degrees of freedom $t,\dot t,x,\dot x$ may have an effect on both the expressions we obtain for the equations of motion and the Noether charges. To investigate whether this is the case, we absorb the Lagrange multiplier terms into the integral of the action using the Dirac delta functions $\delta(\gamma-\gamma_i)$ for the initial conditions and  $\delta(\gamma-\gamma_f)$ for the connecting conditions. This allows us to define two new terms $\tilde E^{(1,2)}_{\rm BVP}$ for the forward and backward path, in which we absorb any reference to $t_1,x_1$ and $t_2,x_2$ including the those related to the Lagrange multiplier terms
\begin{align}
 \nonumber   {\cal E}_{\rm IVP}^{\lambda}&=\int_{\gamma_i}^{\gamma_f} d\gamma \Big\{ E_{\rm BVP}[t_1,\dot t_1,x_1,\dot x_1] - E_{\rm BVP}[t_2,\dot t_2,x_2,\dot x_2] \\
 + \Big(& \nonumber\lambda_1\big({ t}_1(\gamma)-t_i\big)+\lambda_2\big(\dot { t}_1(\gamma)-\dot t_i\big) \\
\nonumber  &\hspace{2cm} +\lambda_3\big({ x}_1(\gamma)-x_i\big) + \lambda_4\big(\dot{ x}_1(\gamma)-\dot x_i\big) \Big)\delta(\gamma-\gamma_i) \\
\nonumber+\Big( &\lambda_5\big({ t}_1(\gamma)-{ t}_2(\gamma)\big) +\lambda_6\big( \dot{ t}_1(\gamma)- \dot{ t}_2(\gamma)\big)\\
 \nonumber &\hspace{2cm}+ \lambda_7\big({ x}_1(\gamma)-{ x}_2(\gamma)\big)+\lambda_8\big( \dot{ x}_1(\gamma)- \dot{ x}_2(\gamma)\big) \Big) \delta(\gamma-\gamma_f)\Big\}\\
&=\int_{\gamma_i}^{\gamma_f} d\gamma \Big\{ \tilde E_{\rm BVP}^{(1)}[t_1,\dot t_1,x_1,\dot x_1] - \tilde E^{(2)}_{\rm BVP}[t_2,\dot t_2,x_2,\dot x_2] \Big\}
 \label{eq:EIVPcontExpl2}.
\end{align}

When deriving the geodesic equations, which are the Euler-Lagrange equations for this action, we need to vary ${\cal E}_{\rm IVP}^{\lambda}$ with respect to all dynamical degrees of freedom, including the Lagrange multipliers. The derivation is most transparent when one goes over to central and relative coordinates $x_+=\frac{1}{2}(x_1+x_2)$ and $x_-=x_1-x_2$ and $t_+=\frac{1}{2}(t_1+t_2)$ and $t_-=t_1-t_2$, so that the action may be written in analogy with \cref{eq:EIVPcontExpl2} as
\begin{align}
{\cal E}_{\rm IVP}^{\lambda}=\int_{\gamma_i}^{\gamma_f} d\gamma \Big\{ \tilde E_{\rm BVP}^{(-)}[t_-,\dot t_-,x_-,\dot x_-] - \tilde E^{(+)}_{\rm BVP}[t_+,\dot t_+,x_+,\dot x_+] \Big\}.
\end{align}
Ref.~\cite{galley_classical_2013,Rothkopf:2022zfb} discuss in detail that the classical solution emerges from the Euler-Lagrange equation in the relative coordinates in the so-called physical limit, where the relative coordinates vanish. The corresponding contribution to the action, including the Lagrange multiplier terms reads
\begin{align}
\nonumber & \tilde E_{\rm BVP}^{(-)}[t_-,\dot t_-,x_-,\dot x_-]=E_{\rm BVP}[t_-,\dot t_-,x_-,\dot x_-] \\
&+ \frac{1}{2} \Big( \nonumber\lambda_1\big({ t}_-(\gamma)-t_i\big)+\lambda_2\big(\dot { t}_-(\gamma)-\dot t_i\big) \\
&\hspace{2.6cm} +\lambda_3\big({ x}_-(\gamma)-x_i\big) + \lambda_4\big(\dot{ x}_-(\gamma)-\dot x_i\big) \Big)\delta(\gamma-\gamma_i) \\
&+\Big( \lambda_5\big({ t}_-(\gamma)\big) +\lambda_6\big( \dot{ t}_-(\gamma)\big) + \lambda_7\big({ x}_-(\gamma)\big)+\lambda_8\big( \dot{ x}_-(\gamma)\big) \Big) \delta(\gamma-\gamma_f)
\end{align}
When varying $\tilde E_{\rm BVP}^{(-)}[t_-,\dot t_-,x_-,\dot x_-]$ the Lagrange multiplier terms lead to contributions including $\delta \dot t_-$ and $\delta \dot x_-$. Carrying out integration-by-parts as required to make reference only to the variation of $\delta t_-$ and $\delta x_-$, the resulting terms involving $\lambda_2,\lambda_6,\lambda_4$ and $\lambda_8$ all contain a $\gamma$ derivative acting on the Delta function. On the other hand, the terms featuring $\lambda_1,\lambda_5,\lambda_3$ and $\lambda_7$ lead to variations with $\delta t_-$ and $\delta x_-$ directly so that no integration by parts is necessary.

The various contributions featuring $\lambda_i$'s, derived above, were not explicitly considered in \cite{Rothkopf:2023ljz} and, as will be shown below, are important to understanding the conservation of the Noether charge. The relevant Euler Lagrange equations for the relative coordinates $t_-$ and $x_-$ are
\begin{align}
&\frac{\partial E_{\rm BVP}}{\partial t_-} - \frac{d}{d\gamma} \frac{\partial E_{\rm BVP}}{\partial \dot t_-} +\frac{1}{2}\lambda_1\delta(\gamma-\gamma_i) +\lambda_5\delta(\gamma-\gamma_f)\\
\nonumber &\hspace{3.2cm}-  \frac{d}{d\gamma}\Big( \frac{1}{2}\lambda_2\delta (\gamma-\gamma_i) + \lambda_6\delta(\gamma-\gamma_f)\Big)=0,\\
&\frac{\partial E_{\rm BVP}}{\partial x_-} -\frac{d}{d\gamma} \frac{\partial E_{\rm BVP}}{\partial \dot x_-} +\frac{1}{2}\lambda_3\delta(\gamma-\gamma_i) +\lambda_7\delta(\gamma-\gamma_f) \\
\nonumber&\hspace{3.2cm}-  \frac{d}{d\gamma}\Big( \frac{1}{2}\lambda_4\delta (\gamma-\gamma_i) + \lambda_8\delta(\gamma-\gamma_f)\Big)=0,
\end{align}
under the constrains enforced by the variation of the Lagrange multipliers $\delta \lambda_i$. We see that in contrast to the implicit treatment of the boundary and connecting conditions, additional terms appear (the one's containing the delta functions and derivatives of delta functions) in the equations of motion that affect the boundary region. The fact that there appears a derivative acting on the delta function at $\gamma_f$ will lead to modifications not just exactly on the boundary but also further inside the simulation domain. The assignment of the Lagrange multiplier terms introduced above however is not unique, due to the fact that $\delta t_1(\gamma_f)=\delta t_2(\gamma_f)$ and $\delta x_1(\gamma_f)=\delta x_2(\gamma_f)$. I.e. all the terms located at $\gamma_f$ can be moved between the forward and backward branch, and consequently the sign in front of the terms $\lambda_{5-8}$ is not uniquely fixed. In addition it is also possible to absorb e.g. the terms with $\lambda_5$ all into one of the branches, where they then identically cancel. I.e. at this stage, there remains an ambiguity about which of the Lagrange multiplier terms affect the boundary region in the geodesic equations. We will identify the relevant terms through an inspection of the actual solution.

The situation is more transparent for the Noether charge corresponding to time translations. If we vary the action ${\cal E}_{\rm IVP}^{\lambda}$ such that the variation of the individual time degrees of freedom is compatible with the forward-backward construction needed for IVPs, i.e. $\delta t_1=\delta t_2$ with $t_1\to t_1+dt$ and $t_2\to t_2+dt$, there is no change in the value of the action $\delta {\cal E}_{\rm IVP}^{\lambda}=0$. This is a consequence of the manifest time translation invariance of the action we constructed. When deriving how the time translations affect the action formally we end up with the following expression 
\begin{align}
\nonumber 0=\delta {\cal E}_{\rm IVP}^{\lambda}= \int d\gamma\Big\{ &\Big(\underbracket{\frac{\partial \tilde E^{(1)}_{\rm BVP}}{\partial t_1} - \partial_\gamma \frac{\partial \tilde E^{(1)}_{\rm BVP}}{\partial \dot t_1}}_{\rm e.o.m. \& \lambda\,{\rm terms}=0}\Big)\delta t_1 + \partial_\gamma \Big(\underbracket{\frac{\partial \tilde E^{(1)}_{\rm BVP}}{\partial \dot t_1}}_{Q_{t,1}^{\lambda}} \delta t_1 \Big)\\
- &\Big(\underbracket{\frac{\partial \tilde E^{(2)}_{\rm BVP}}{\partial t_2} - \partial_\gamma \frac{\partial \tilde E^{(2)}_{\rm BVP}}{\partial \dot t_2}}_{\rm e.o.m. \& \lambda\,{\rm terms}=0}\Big)\delta t_2 - \partial_\gamma \Big( \underbracket{\frac{\partial \tilde E^{(2)}_{\rm BVP}}{\partial \dot t_2}}_{Q_{t,2}^{\lambda}} \delta t_2\Big)
 \Big\}.\label{eq:Noether} \end{align}
Since $\delta t_1=dt=\delta t_2$, the Euler-Lagrange like terms in the first and second row cancel against each other, including the Lagrange multiplier terms. This leaves only the two total derivates that encode the Noether charge on the forward $Q_{t,1}^{\lambda}$ and backward branches $Q_{t,2}^{\lambda}$. Note that in the forward-backward construction \cref{eq:Noether} only states that the difference between the Noether charge on the forward and backward branch vanishes.

We emphasize here that in the presence of the explicit $\lambda_i$ terms it is not  $ E_{\rm BVP}$ but $\tilde E_{\rm BVP}$ which appears in the definition of the Noether charges $Q_t^\lambda$. Writing out the terms explicitly we find that the expression for the Noether charge is modified with respect to the case, where the connecting conditions were implicitly enforced. Since the definition of the Noether charge on the forward and backward branch asks us to take the derivative with respect to $\dot t_{1,2}$ respectively, only the terms involving $\lambda_2$ and $\lambda_6$ can contribute. The explicit expression for the Noether charge on forward branch according to \cref{eq:Noether} is
\begin{align}
Q_{t,1}^{\lambda} =  \frac{\partial \tilde E^{(1)}_{\rm BVP}}{\partial \dot t_{1}} = \frac{\partial E_{\rm BVP}}{\partial \dot t_{1}} +  \lambda_2  \delta (\gamma-\gamma_i) +\lambda_6  \delta (\gamma-\gamma_f).\label{eq:modNoether}
\end{align}
We will demonstrate that it is indeed \cref{eq:modNoether} which is conserved over the whole world-line in the following section.

 \section{Discrete formalism and numerical results}
 
To discretize the action functional ${\cal E}_{\rm IVP}^\lambda$  along the world-line parameter $\gamma\in[\gamma_i,\gamma_f]$ with $N_\gamma$ steps, we use an equidistant $d\gamma=(\gamma_f-\gamma_i)/(N_\gamma-1)$. The discretized forward and backward paths $x_{1},x_{2}$ and times $t_{1},t_{2}$ are referred to as ${\bf x}_{1,2}=(x_{1,2}(0),x_{1,2}(\Delta \gamma),x_{1,2}(2\Delta \gamma),\ldots,x_{1,2}((N_\gamma-1)\Delta\gamma))^{\rm T}$ and ${\bf t}_{1,2}=(t_{1,2}(0),t_{1,2}(\Delta \gamma),t_{1,2}(2\Delta \gamma),\ldots,t_{1,2}((N_\gamma-1)\Delta\gamma))^{\rm T}$ respectively.

Deploying appropriately regularized summation-by-parts finite difference operators ${\bar{\mathds{D}}^{\rm R}}$ (here we will use second order \texttt{SBP21} and fourth order \texttt{SBP42} SBP operators) together with a compatible quadrature rule  described by the matrix ${\bar{\mathds{H}}}$ (for details of the construction of these operators see \cite{Rothkopf:2022zfb,Rothkopf:2023ljz})  we can write down the discretized counterpart of \cref{eq:EIVPcontExpl2} as
\begin{align}
\nonumber \mathds{E}_{\rm IVP}= &\frac{1}{2} \left\{ ({\bar{\mathds{D}}^{\rm R}}_t{\bf t}_1)^{\rm T} \mathbb{d}\left[c^2+\frac{2 {\bf V}({\bf x}_1)}{m}\right] {\bar{\mathds{H}}} ({\bar{\mathds{D}}^{\rm R}}_t{\bf t}_1) -  ({\bar{\mathds{D}}^{\rm R}}_x{\bf x}_1)^{\rm T} {\bar{\mathds{H}}} ({\bar{\mathds{D}}^{\rm R}}_x{\bf x}_1)\right\}\\
\nonumber-&\frac{1}{2} \left\{({\bar{\mathds{D}}^{\rm R}}_t{\bf t}_2)^{\rm T} \mathbb{d}\left[c^2+\frac{2 {\bf V}({\bf x}_2)}{m}\right] {\bar{\mathds{H}}} ({\bar{\mathds{D}}^{\rm R}}_t{\bf t}_2) -  ({\bar{\mathds{D}}^{\rm R}}_x{\bf x}_2)^{\rm T} {\bar{\mathds{H}}} ({\bar{\mathds{D}}^{\rm R}}_x{\bf x}_2)\right\}\\
\nonumber+&\lambda_1\big({\bf t}_1[1]-t_i\big)+\lambda_2\big((\mathds{D}{\bf t}_1)[1]-\dot t_i\big)+\lambda_3\big({\bf x}_1[1]-x_i\big)+ \lambda_4\big((\mathds{D}{\bf x}_1)[1]-\dot x_i\big)\\
\nonumber +&\lambda_5\big({\bf t}_1[N_\gamma]-{\bf t}_2[N_\gamma]\big) + \lambda_6\big( (\mathds{D}{\bf t}_1)[N_\gamma]- (\mathds{D}{\bf t}_2)[N_\gamma]\big)\\
+&\lambda_7\big({\bf x}_1[N_\gamma]-{\bf x}_2[N_\gamma]\big)+\lambda_8\big( (\mathds{D}{\bf x}_1)[N_\gamma]- (\mathds{D}{\bf x}_2)[N_\gamma]\big)\label{eq:discEIVP}.
\end{align}
The above expression implies matrix vector multiplication, whenever a matrix quantity such as $\bar{\mathds{H}}$ or $\bar{\mathds{D}}$ acts on a vector ${\bf x}_{1,2}$ or ${\bf t}_{1,2}$. The matrix expression $\mathbb{d}[f({\bf x})]$ contains on the diagonal the values $\mathbb{d}_{kk}=f({\bf x}(\gamma_k))$ and zero otherwise.
 
When incorporating the initial and connecting conditions of the forward-backward path construction explicitly, we showed in the previous section that the Noether charge receives correction terms proportional to a delta function at the boundaries of the $\gamma$ domain. We must therefore define the discrete counterpart of the Dirac delta function $ \mathfrak{d}_k$ (also called lifting operator when applied at the boundary \cite{nordstrom2017roadmap}). For an arbitrary discretized function ${\bf f}$, when integrated over $ \mathfrak{d}[k]$ we must obtain the value of the function at the $k$-th $\gamma$ step, i.e.   ${\bf f}{\mathds{H}}\mathfrak{d}[k] ={\bf f}[k]$. This motivates the discrete definition $ \mathfrak{d}[k] = {\mathds{H}}^{-1} {\bf e}_k $, with ${\bf e}_k$ the unit vector with entry one at position $k$.
 
  \subsection{Quartic potential example}
 
Choosing as concrete potential the function $V(x)=\kappa x^4$, we discretize the world-line of the particle motion between $\gamma_i=0$ and $\gamma_f=1$ with $N_\gamma=32$ points. Without loss of generality, we arbitrarily set the starting time to $t_i=0$ and the starting position to $x_i=1$. To obtain an initial velocity $v_i=1/10$ we choose $\dot t=1$ and $\dot x=v_i$. Note that we do not fix the value of $t_f$ but only the initial velocity of time with respect to $\gamma$. The choice of $\dot t=1$ will lead to dynamics, such that $t_f$ will be of the order of one. As strength for the linear potential we choose $\alpha=1/4$. The corresponding discrete action functional reads explicitly
\begin{align}
\nonumber \mathds{E}^{\rm qrt}_{\rm IVP}= &\frac{1}{2} \left\{ ({\bar{\mathds{D}}^{\rm R}}_t{\bf t}_1)^{\rm T}\mathbb{d}\left[1+2 \kappa {\bf x}^4_1\right] {\bar{\mathds{H}}} ({\bar{\mathds{D}}^{\rm R}}_t{\bf t}_1) -  ({\bar{\mathds{D}}^{\rm R}}_x{\bf x}_1)^{\rm T} {\bar{\mathds{H}}} ({\bar{\mathds{D}}^{\rm R}}_x{\bf x}_1)\right\}\\
\nonumber-&\frac{1}{2} \left\{ ({\bar{\mathds{D}}^{\rm R}}_t{\bf t}_2)^{\rm T} \mathbb{d}\left[1+2 \kappa {\bf x}^4_2\right] {\bar{\mathds{H}}} ({\bar{\mathds{D}}^{\rm R}}_t{\bf t}_2) -  ({\bar{\mathds{D}}^{\rm R}}_x{\bf x}_2)^{\rm T} {\bar{\mathds{H}}} ({\bar{\mathds{D}}^{\rm R}}_x{\bf x}_2)\right\}\\
\nonumber+&\lambda_1\big({\bf t}_1[1]-t_i\big)+\lambda_2\big((\mathds{D}{\bf t}_1)[1]-\dot t_i\big) +\lambda_3\big({\bf x}_1[1]-x_i\big)+\lambda_4\big((\mathds{D}{\bf x}_1)[1]-\dot x_i\big)\\
\nonumber+&\lambda_5\big({\bf t}_1[N_\gamma]-{\bf t}_2[N_\gamma]\big) +\lambda_6\big( (\mathds{D}{\bf t}_1)[N_\gamma]- (\mathds{D}{\bf t}_2)[N_\gamma]\big)\\
 +&\lambda_7\big({\bf x}_1[N_\gamma]-{\bf x}_2[N_\gamma]\big)+\lambda_8\big( (\mathds{D}{\bf x}_1)[N_\gamma]- (\mathds{D}{\bf x}_2)[N_\gamma]\big) \label{eq:discEIVPqrt},
\end{align}
where taking the fourth power of the ${\bf x}_{1,2}$ vector is to be understood in an element wise fashion.
 
In \cite{Rothkopf:2023ljz} we discussed in detail the solutions $t(\gamma)$ and $x(\gamma)$, as well as the physical trajectory $x(t)$, one obtains from numerically determining the critical point of \cref{eq:discEIVPqrt}. We found that even though the $\gamma$ grid is equidistantly spaced, a non-equidistant spacing for the time coordinate emerges dynamically. It adapts to the dynamics of the spatial coordinate, showing a finer temporal resolution, where more curvature exists in $x$. This constitutes a form of automatic adaptive mesh refinement (AMR). 
 
\begin{figure}
\centering
    \hspace{-0.8cm}\includegraphics[scale=0.32]{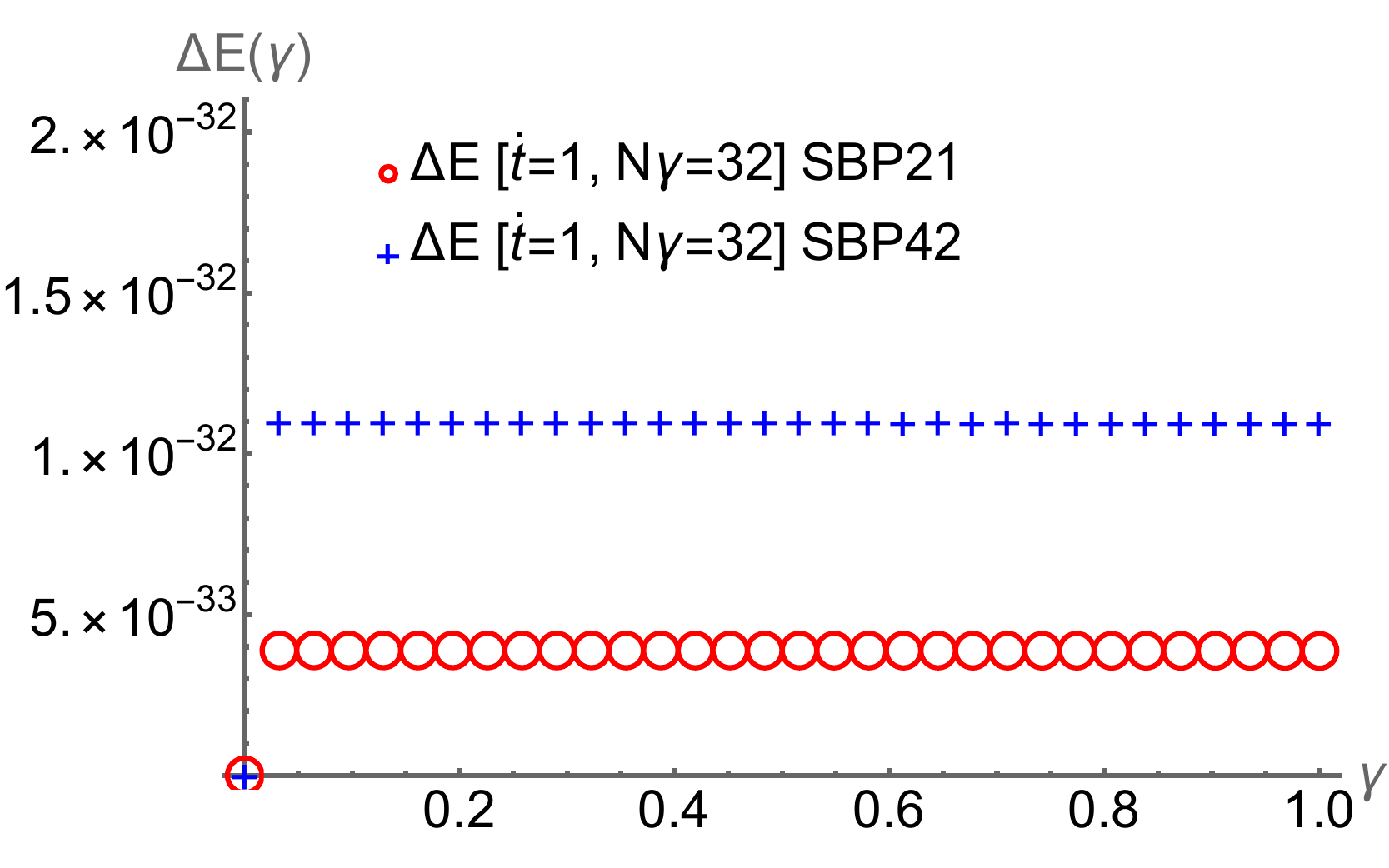}
     \includegraphics[scale=0.3]{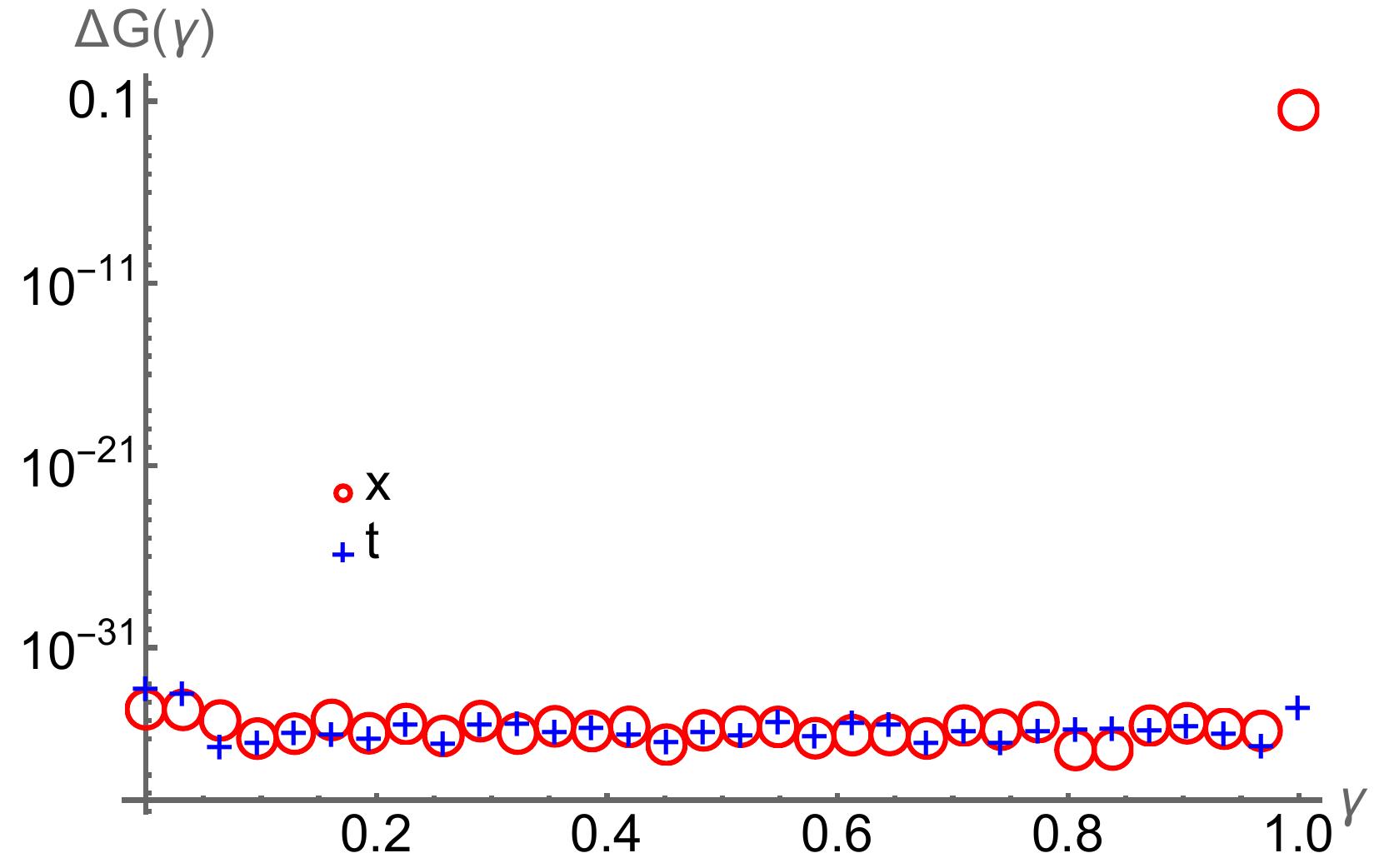}
    \caption{ (left) Deviation $\Delta {\bf E}$ from its continuum value, as obtained from the critical point of $\mathds{E}^{\rm qrt}_{\rm IVP}$ with $V(x)=\kappa x^4$. Results for the \texttt{SBP21} operator are given as red circles and those for \texttt{SBP42} as blue crosses. Note that the deviation $\Delta{\bf E}$ is exactly zero within machine precision. (right) Deviation of the spatial coordinate $\Delta{\bf G}^x$ (red circles) and time coordinate $\Delta{\bf G}^t$ (blue crosses) from the discretized geodesic equation.}
    \label{fig:NLDeltaE}
\end{figure}

 Here we focus on the central quantity of interest in this study, $Q_t^{\lambda}$, defined in \cref{eq:modNoether}, which in the continuum represents the conserved charge associated with the time-translation symmetry of the system. We consider it here in its naively discretized form\begin{align}
    {\bf \mathds{Q}_t}=(\mathds{D}{\bf t})\circ({\bf 1} + 2\kappa {\bf x}^4) +\lambda_2\mathfrak{d}_1 +\lambda_6 \mathfrak{d}_{N_\gamma} ,
\end{align}
and we define its deviation from the continuum result via the difference
\begin{align}
    \Delta {\bf E} = {\bf \mathds{Q}_t} -Q_t = (\mathds{D}{\bf t})\circ({\bf 1} + 2 \kappa {\bf x}^4) +\lambda_2\mathfrak{d}_1 +\lambda_6 \mathfrak{d}_{N_\gamma} - \dot t_i (1+2\kappa x_i^4),
\end{align}
which we plot in \cref{fig:NLDeltaE} using the \texttt{SBP21} operator (red circles) and the \texttt{SBP42} operator (blue crosses).

First and foremost, the discretized Noether charge ${\bf \mathds{Q}_t}^\lambda$ is \textit{exactly conserved} in the discrete setting over the whole of the simulated time domain. Since we use a minimizer to find the critical point with \texttt{WorkingPrecision} set to 40, the values of $<10^{-30}$ indeed reflect a true zero.

Secondly, the value of ${\bf \mathds{Q}_t}^\lambda$ that remains conserved in the interior agrees with the true continuum value, prescribed by the initial conditions, \textit{within machine precision}. This is a highly non-trivial result, as even in energy preserving schemes, such as the leap-frog, the conserved quantities do not necessarily agree with the continuum ones. We find that the Lagrange multiplier $\lambda_2$ at the critical point of $\mathds{E}^{\rm qrt}_{\rm IVP}$ becomes zero. In turn ${\bf \mathds{Q}_t}^\lambda$ takes on the continuum value by construction at the first point in $\gamma$, as there it is defined solely by the provided initial conditions. Since its value is conserved over the whole domain, it actually exhibits the correct continuum value everywhere.
 
In \cite{Rothkopf:2023ljz} we found that ${\bf x}$ and ${\bf t}$, obtained from the critical point of \cref{eq:discEIVPqrt} fulfil the discretized Euler-Lagrange equations associated with \cref{eq:discEIVPqrt} in the interior of the $\gamma$ interval, except for the last two points. We discussed above that we understand this deviation close to the boundary to arise from corrections introduced by the Lagrange multipliers. However we also found in \cref{sec:expllagrange} that there exists an ambiguity in assigning Lagrange multiplier contributions to the forward and backward branches due to the fact that variations at the end of the paths agree.

 We attempted to infer the appropriate boundary corrections by trial and error and succeeded for the time coordinate, which obeys the following discrete geodesic equation
\begin{align}
&  \mathds{D}\big( (1+2\kappa {\bf x}^4 )\circ\mathds{D}{\bf t} \big) + \lambda_6\mathds{D}\mathfrak{d}_{N_\gamma}=\Delta {\bf G}^t. \label{eq:DiscrGeodtLag}
\end{align}
When inserted into this equation, the solution ${\bf t}$, found at the critical point of \cref{eq:discEIVPqrt} produces $\Delta {\bf G}^t=0$ everywhere in the whole simulation domain as shown by the blue crosses in the right panel of \cref{fig:NLDeltaE}. 

The situation is not completely satisfactory for the spatial coordinate, where we were able to account for the correction to the second to last point in $\gamma$ but were unable to find the expression for the correction at the very last point
\begin{align}
&\mathds{D}\mathds{D}{\bf x} + (4\kappa {\bf x}^3) \circ (\mathds{D}{\bf t})\circ (\mathds{D}{\bf t}) - \lambda_8\mathds{D}\mathfrak{d}_{N_\gamma}=\Delta {\bf G}^x. \label{eq:DiscrGeodxLag}
\end{align}
As shown by the red circles in the right panel of \cref{fig:NLDeltaE}, inserting the solution ${\bf x}$ into the above equation, we find that $\Delta {\bf G}^x$ vanishes within the interior of the simulation domain but a small correction remains at the final $\gamma$ step. 

The exact zero values of $\Delta {\bf G}^t$ and $\Delta {\bf G}^x$ in the interior, shown in the right panel of \cref{fig:NLDeltaE}, confirm that there our solution obtained directly from the action principle fulfils the naively discretized geodesic equations and only at the boundary corrections receives modifications from the Lagrange multipliers. We showed in  \cite{Rothkopf:2023ljz} that these boundary corrections do not spoil the favorable convergence to the continuum solution under grid refinement.

Since the form of the boundary contributions does not depend on the potential used in the action, the same boundary corrections as in the quartic potential case apply also for the linear potential discussed in \cite{Rothkopf:2023ljz}. We have explicitly checked that the Noether charge for the linear potential too is preserved in the whole simulation domain.
 
 \section{Summary and conclusion}
 
We have shown that the Langrange multipliers, used in a recently developed variational discretization approach to IVPs, lead to correction terms affecting the equations of motion and the Noether charges of the system. These contributions were not considered in the previous publication \cite{Rothkopf:2023ljz}. By taking them into account, we are able to derive an expression for the Noether charge that is exactly preserved at the true continuum value in the whole simulation domain. After investigating possible boundary terms to the naively discretized Euler-Lagrange equations, we aposteriori selected a subset of these, which lead to modified Euler-Lagrange equations. The classical trajectory ${\bf t}$ indeed fulfils these equations of motions everywhere and only ${\bf x}$ shows a single remaining deviating point at the final $\gamma$ step. This short letter thus provides a clarification and improvement over the observed deviations in \cite{Rothkopf:2023ljz}.
  
 \section*{Acknowledgements}
 A.R. thanks Will Horowitz for insightful discussions. J.~N. was supported by the Swedish Research Council grant nr. 2021-05484. 
 
\FloatBarrier

\begin{backmatter}

\section*{Competing interests}
  The authors declare that they have no competing interests.

\section*{Author's contributions}
    \begin{itemize}
         \item A. Rothkopf: derivation of the correction terms for the Noether charge and heuristic exploration of correction terms for the equations of motion
         \item J. Nordst\"om: guidance on the formulation and implementation of SBP based discretization schemes, literature review, editing
    \end{itemize}


\bibliographystyle{stavanger-mathphys}


\bibliography{refs}


\end{backmatter}


\end{document}